\documentclass[conference]{IEEEtran}
\IEEEoverridecommandlockouts

\usepackage[utf8]{inputenc}
\usepackage[T1]{fontenc}

\usepackage{amsmath,amsfonts,amssymb,amsthm}
\usepackage{mathtools}

\usepackage{enumerate}
\usepackage{cite}
\usepackage{url}
\usepackage{xcolor}

\usepackage{graphicx}
\DeclareGraphicsExtensions{.eps}
\usepackage{caption}
\usepackage{subcaption}

\usepackage{calc}

\usepackage{tikz}
\usetikzlibrary{shapes, decorations.markings, decorations.pathreplacing, calligraphy}

\newtheorem{Theorem}{Theorem}
\newtheorem{Proposition}[Theorem]{Proposition}

\newtheorem{Assumption}{Assumption}

\theoremstyle{remark}
\newtheorem{Remark}{Remark}

\def \ba{\begin{array}}
\def \ea{\end{array}}
\def \be{\begin{equation}}
\def \ee{\end{equation}}

\def \colsep{\arraycolsep}

\def \bb{\mathbb}

\def \mc{\mathcal}

\def \ms{\mathsf}

\def \diag{\mathrm{diag}}

\def \eig{\mathrm{eig}}

\def \cond{\mathrm{cond}}

\def \test{\mathrm{test}}
\def \train{\mathrm{train}}
\def \t{\mathtt{T}}

\DeclareMathOperator*{\esssup}{ess\,sup}

\definecolor{true_in}{rgb}{0.6350 0.0780 0.1840}
\definecolor{est_in}{rgb}{0.4660 0.6740 0.1880}
\definecolor{true_out}{rgb}{0.9290 0.6940 0.1250}
\definecolor{est_out}{rgb}{0 0.4470 0.7410}
\definecolor{grey}{RGB}{211,211,211}

\usepackage{geometry}
\newgeometry{top=1in,bottom=0.75in,right=0.75in,left=0.75in}

\begin{document}

\title{Learning-based Design of Luenberger Observers for Autonomous Nonlinear Systems}

\author{Muhammad Umar B. Niazi$^{*,\dagger}$ \and 
John Cao$^*$ \and
Xudong Sun$^*$ \and 
Amritam Das$^*$ \and 
Karl Henrik Johansson$^*$
\thanks{$*$ Division of Decision and Control Systems and Digital Futures, EECS, KTH Royal Institute of Technology, SE-100 44 Stockholm, Sweden.}
\thanks{$\dagger$ Laboratory for Information and Decision Systems, Massachusetts Institute of Technology, 77 Massachusetts Avenue, Cambridge, MA 02139, USA. Corresponding author's email: \texttt{niazi@mit.edu}}
\thanks{This work is supported by the Swedish Research Council and the Knut and Alice Wallenberg Foundation, Sweden. It also received funding from European Union's Horizon Research and Innovation Programme under Marie Sk\l{}odowska-Curie grant agreement No. 101062523.}
}

\maketitle
\thispagestyle{empty}
\pagestyle{empty}

\begin{abstract}
Designing Luenberger observers for nonlinear systems involves the challenging task of transforming the state to an alternate coordinate system, possibly of higher dimensions, where the system is asymptotically stable and linear up to output injection. The observer then estimates the system's state in the original coordinates by inverting the transformation map. However, finding a suitable injective transformation whose inverse can be derived remains a primary challenge for general nonlinear systems. We propose a novel approach that uses supervised physics-informed neural networks to approximate both the transformation and its inverse. Our method exhibits superior generalization capabilities to contemporary methods and demonstrates robustness to both neural network's approximation errors and system uncertainties.
\end{abstract}


\section{Introduction}

Nonlinear Luenberger observers, also known as Kazantzis-Kravaris/Luenberger (KKL) observers, generalize the theory of Luenberger observers \cite{luenberger1964} to nonlinear systems. The main idea of KKL observers is to find an injective map that satisfies a certain partial differential equation (PDE) and transforms a nonlinear system to another coordinate system, possibly of higher dimensions than the original state space. The dynamics of the transformed system are required to be stable and linear up to output injection. Then, the KKL observer is a copy of the transformed system and estimates the state of the original system by inverting the transformation map.

Initially proposed by \cite{shoshitaishvili1990} and \cite{shoshitaishvili1992}, the theory of KKL observers was subsequently rediscovered by Kazantzis \& Kravaris \cite{kazantzis1998}, who provided local guarantees around an equilibrium point via Lyapunov's Auxiliary Theorem. Although \cite{krener2002b} relaxed the restrictive assumptions of \cite{kazantzis1998} to some extent, the analysis remained local until \cite{kreisselmeier2003} proposed the first global result under the assumption of the so-called finite complexity, which also turned out to be quite restrictive for general nonlinear systems. In this regard, a complete and most general treatment of the problem was presented by Andrieu \& Praly \cite{andrieu2006}, who introduced the notion of \textit{backward distinguishability} ensuring the existence of an injective transformation required by the KKL observers. Later, under some additional observability conditions, \cite{andrieu2014} proved that KKL observers converge exponentially and are also tunable. The theory is also extended to non-autonomous and controlled nonlinear systems in \cite{engel2007, bernard2017, bernard2018}.

The main challenge in the design of KKL observers is to not only find the transformation map but also its left inverse, and both problems turn out to be very difficult in practice; see \cite{andrieu2021} and \cite{bernard2022}. To this end, \cite{ramos2020, peralez2021, buisson2022} have proposed several methods to approximate the transformation map and its inverse via feedforward neural networks. By fixing the dynamics of the KKL observer, they propose to generate synthetic data trajectories by numerically solving both the system's model and the KKL observer, where both are initialized at multiple points in their corresponding state spaces. Then, using a supervised learning approach, a neural network is trained to approximate the transformation map and its left inverse. Moreover, \cite{buisson2022} also proposed an unsupervised learning approach by assuming an autoencoder-type architecture and adding the PDE associated with the transformation map as a design constraint. However, both approaches suffer from overfitting on the training samples and do not generalize well in practice.

In this paper, we propose a \emph{supervised physics-informed learning} approach to approximate the transformation map and its left inverse. Such an approach incorporates the physical knowledge described by the PDE constraint, which is directly integrated with the conventional supervised learning \cite{raissi2019, karniadakis2021}. Embedding the physical knowledge of systems by adding the PDE constraint as a physically relevant invariant improves the accuracy, generalization, and training time of the learning method. In this way, we improve upon the idea of \cite{ramos2020, peralez2021, buisson2022} by avoiding overfitting and obtaining better generalization to the whole state space.

The main contribution of this paper includes a complete learning method of the KKL observer design via a supervised physics-informed neural network (PINN). We show that the KKL observer is robust to not only the neural network's approximation error but also to model and sensor uncertainties. The robustness is quantified in terms of input-to-state stability \cite{sontag1995} of the state estimation error. We define an empirical metric to quantify the generalization capability of the learned KKL observer and provide a detailed discussion on why our method exhibits better generalization capabilities than the supervised neural network (NN) approach of \cite{ramos2020,peralez2021,buisson2022} and the unsupervised autoencoder (AE) approach of \cite{buisson2022}. Finally, we demonstrate the dominance of our method over these approaches through statistically well-designed experiments.

After summarizing a general idea of KKL observers in Section~\ref{sec_KKLprelim}, we state the problem addressed in this paper in Section~\ref{sec_prob}. The learning method of KKL observers is presented in Section~\ref{sec_learning}. Section~\ref{sec_eval} evaluates the performance of the observer under approximation errors and uncertainties, and defines and discusses an empirical metric to assess the generalization capability of the learned observer. Finally, Section~\ref{sec_exp} presents the experimental results and Section~\ref{sec_conc} ends with concluding remarks and the future outlook.

\textit{Notations.} For a vector $x\in\bb{R}^n$, the Euclidean norm $\|x\|=\sqrt{x^\t x}$ and the maximum norm $\|x\|_\infty=\max_i |x_i|$.
For a measurable essentially bounded function $w\in L^\infty(\bb{R};\bb{R}^n)$, the essential supremum norm $\|w\|_{L^\infty}=\esssup_{t\in\bb{R}} \|w(t)\| \doteq  \inf\{c\geq 0: \|w(t)\|_\infty\leq c ~\text{for almost every}~ t\in\bb{R}\}$. 
For a matrix $M\in\bb{R}^{n\times m}$, $\|M\|=\sup_{\|x\|=1} \|Mx\|$ denotes the induced norm, which is equal to the maximum singular value $\sigma_{\max}(M)$. The spectrum of $M\in\bb{R}^{n\times n}$ is denoted by $\eig(M)$, and $\lambda_{\min}(M) = \min_{\lambda\in\eig(M)} |\text{Re}(\lambda)|$ and $\lambda_{\max}(M) = \max_{\lambda\in\eig(M)} |\text{Re}(\lambda)|$. The condition number of $M$ is denoted by $\cond(M)$.

\section{Preliminaries on KKL Observers}
\label{sec_KKLprelim}

In this section, we briefly summarize the theory of KKL observers. For more details, see \cite{kazantzis1998, andrieu2006, andrieu2014}.

Consider a nonlinear system
\be \label{eq:sys_x}
    \dot{x} = f(x); \quad y = h(x)
\ee
where $x(t)\in\mc{X}\subset\bb{R}^{n_x}$ is the state with $x(0)=x_0\in\mc{X}$ the initial condition, $y(t)\in\bb{R}^{n_y}$ is the measured output, and the maps $f:\mc{X}\rightarrow\bb{R}^{n_x}$ and $h:\mc{X}\rightarrow\bb{R}^{n_y}$ are smooth.

The design method of a KKL observer is as follows:
\begin{enumerate}
\item Find an injective\footnote{A map $\mc{T}:\mc{X}\rightarrow\bb{R}^{n_z}$ is said to be \textit{injective} if for every $x_1,x_2\in\mc{X}$, $\mc{T}(x_1)=\mc{T}(x_2)$ implies $x_1=x_2$.} 
map $\mc{T}:\mc{X}\rightarrow\bb{R}^{n_z}$ that transforms \eqref{eq:sys_x} to new coordinates $z=\mc{T}(x)$, where
\be \label{eq:sys_z}
    \dot{z} = Az + Bh(x); \quad z(0)=\mc{T}(x_0)
\ee
with $A\in\bb{R}^{n_z\times n_z}$ a Hurwitz matrix and $B\in\bb{R}^{n_z\times n_y}$ such that $(A,B)$ is controllable. 
From \eqref{eq:sys_z}, it follows that $\mc{T}$ must be a solution to the following PDE:
\be \label{eq:pde}
    \frac{\partial \mc{T}}{\partial x}(x) f(x) = A\mc{T}(x) + Bh(x); \quad \mc{T}(0)=0.
\ee
\item Since $\mc{T}$ is injective, its left inverse $\mc{T}^*$ exists, i.e., $\mc{T}^*(\mc{T}(x))=x$. The KKL observer is then given by
\be \label{eq:kkl_observer}
    \colsep=2pt\ba{ccl}
        \dot{\hat{z}} &=& A\hat{z} + By; \quad \hat{z}(0) = \hat{z}_0 \\
        \hat{x} &=& \mc{T}^*(\hat{z}).
    \ea
\ee
\end{enumerate}

There are certain conditions that system \eqref{eq:sys_x} needs to satisfy in order to ensure the existence of a KKL observer \eqref{eq:kkl_observer} in a sense that $\lim_{t\rightarrow\infty}\|\hat{x}(t)-x(t)\|=0$.
Let $x(t;x_0)$ denote the state trajectory of \eqref{eq:sys_x} with $x(0)=x_0$. Then, \eqref{eq:sys_x} is said to be \textit{forward complete within $\mc{X}$} if for every $x_0\in\mc{X}$, $x(t;x_0)\in\mc{X}$ is well-defined for every $t\in\bb{R}_{\geq 0}$. 

\begin{Assumption} \label{assumption_1}
    There exists a compact set $\mc{X}\subset\bb{R}^{n_x}$ such that the system \eqref{eq:sys_x} is forward complete within $\mc{X}$.
\end{Assumption}

A map $\mc{T}:\mc{X}\rightarrow\bb{R}^{n_z}$ is said to be \textit{uniformly injective} if there exists a class $\mc{K}$ function\footnote{A function $\rho:\bb{R}_{\geq 0}\rightarrow\bb{R}_{\geq 0}$ is of class $\mc{K}$ if it is continuous, zero at zero, and strictly increasing. 
} 
$\rho$ such that, for every $x_1,x_2\in\mc{X}$, $\|x_1-x_2\|\leq \rho(\|\mc{T}(x_1)-\mc{T}(x_2)\|)$. 

For the existence of a KKL observer \eqref{eq:kkl_observer}, it is sufficient that \eqref{eq:sys_x} is forward complete and the map $\mc{T}$ satisfying \eqref{eq:pde} is uniformly injective, see \cite[Theorem 1]{andrieu2006}.
Since $A$ is a Hurwitz matrix,
$
\|\hat{z}(t)-z(t)\|=\|\mc{T}(\hat{x}(t))-\mc{T}(x(t))\|
$
converges to zero exponentially. Thus, the uniform injectivity
\be \label{eq:uniformly_injective}
    \|\hat{x}(t)-x(t)\|\leq \rho(\|\mc{T}(\hat{x}(t)) - \mc{T}(x(t))\|)
\ee
implies that $\|\hat{x}(t)-x(t)\|$ also converges to zero. 
However, only asymptotic convergence of the estimation error can be guaranteed because the inverse $\mc{T}^*$ is a nonlinear map, which may destroy the exponentiality of the convergence.

Given an open set $\mc{O}\supset\mc{X}$, the system \eqref{eq:sys_x} is said to be \textit{backward $\mc{O}$-distinguishable} on $\mc{X}$ if for every pair of distinct initial conditions $x_0^1,x_0^2\in\mc{X}$, there exists $\tau<0$ such that $x(t;x_0^1), x(t;x_0^2)\in\mc{O}$ are well-defined for $t\in[\tau,0]$, and
\[
    h(x(\tau;x_0^1))\neq h(x(\tau;x_0^2)).
\]
In other words, this means that there exists a finite negative time such that the output maps, corresponding to different trajectories initialized in $\mc{X}$, can be distinguished before any of the trajectories leaves $\mc{O}$ in backward time.

\begin{Assumption} \label{assumption_2}
    There exists an open bounded set $\mc{O}\supset\mc{X}$ such that \eqref{eq:sys_x} is backward $\mc{O}$-distinguishable on $\mc{X}$.
\end{Assumption}

It turns out that Assumptions~\ref{assumption_1} and \ref{assumption_2} are sufficient for the existence of an injective map $\mc{T}$ satisfying \eqref{eq:pde}. This result is obtained in \cite{andrieu2006, bernard2022, brivadis2022}, which can be restated as follows:

\begin{Theorem}
\label{thm:BOD}
    Let Assumptions~\ref{assumption_1} and \ref{assumption_2} hold. Then, for any controllable $(A,B)\in(\bb{R}^{n_z\times n_z},\bb{R}^{n_z\times n_y})\setminus \mc{J}$ such that $n_z=n_y(2n_x+1)$, $A+\delta I_{n_z}$ is Hurwitz for some $\delta>0$, and $\mc{J}\subset(\bb{R}^{n_z\times n_z},\bb{R}^{n_z\times n_y})$ is a set of zero Lebesgue measure, there exists a uniformly injective map $\mc{T}:\mc{X}\rightarrow\bb{R}^{n_z}$ satisfying \eqref{eq:pde}.
\end{Theorem}

By relying on Theorem~\ref{thm:BOD}, we propose a learning method for $\mc{T}$ and $\mc{T}^*$ under the constraint that $\mc{T}$ satisfies \eqref{eq:pde}. In what follows, we choose and fix $A\in\bb{R}^{n_z\times n_z}$ and $B\in\bb{R}^{n_z\times n_y}$ such that $A$ is Huwitz and $(A,B)$ is controllable, where $n_z = n_y(2n_x+1)$.

\section{Problem Statement}
\label{sec_prob}

We aim to design a KKL observer for \eqref{eq:sys_x} that estimates the state $x(t)$ by using the knowledge of the system's output $y$ and its model $f(\cdot)$ and $h(\cdot)$. That is, the observer \eqref{eq:kkl_observer} ensures 
$
    \lim_{t\rightarrow\infty}\|\hat{x}(t)-x(t)\|=0
$
when $\mc{T}$ and $\mc{T}^*$ are known. In case, $\mc{T}$ and $\mc{T}^*$ are respectively approximated by $\hat{\mc{T}}$ and $\hat{\mc{T}}^*$, then the asymptotic estimation error satisfies
\[
    \limsup_{t\rightarrow\infty}\|\hat{x}(t) - x(t)\| \leq \epsilon
\]
where $\epsilon>0$ depends on the approximation error.

The problem can be divided into two parts:
\begin{enumerate}
    \item Learn the map $\mc{T}$ satisfying the PDE \eqref{eq:pde} and its left inverse $\mc{T}^*$.
    \item Evaluate the performance of the KKL observer in terms of its robustness to the approximation error, model uncertainties, and measurement noise, and its generalization capability when $\mc{T}$ and $\mc{T}^*$ are learned on a discrete subset of $\mc{X}$.
\end{enumerate}

\section{Learning the Transformation Map and its Left Inverse}
\label{sec_learning}

A critical step of KKL observer design is to find the injective map $\mc{T}: \mc X\to \bb{R}^{n_z}$ satisfying the PDE \eqref{eq:pde}, so that \eqref{eq:sys_x} admits a linear representation \eqref{eq:sys_z}, and its left inverse $\mc{T}^*$, so that a state estimate can be obtained in the original state space coordinates. This amounts to solving the PDE \eqref{eq:pde} for $\mc{T}$, whose solution is obtained in \cite{andrieu2006} as
\be \label{eq:T_solution}
\mc{T}(x) = \int_{-\infty}^0 \exp(A\tau) Bh(\breve{x}(\tau;x)) \text{d}\tau
\ee
where $\breve{x}(\tau;x)\in\mc{X}$ is the backward solution initialized at $x\in\mc{X}$, for $\tau\leq 0$, to the modified dynamics $\dot{\breve{x}}(\tau)=g(\breve{x}(\tau))$ with $g(\breve{x}(\tau))=f(\breve{x}(\tau))$ if $\breve{x}(\tau)\in\mc{X}$ and $g(\breve{x}(\tau))=0$ otherwise. However, there are two issues with this solution:
\begin{itemize}
    \item It is practically impossible to obtain a backward output map $h(\breve{x}(\tau;x))$ for $\tau<0$ and then compute the integral \eqref{eq:T_solution} for every initial point $x\in\mc{X}$; \cite{bernard2022}.
    \item Even if $\mc{T}$ is known in any other form\footnote{See \cite{bernard2018} for some of the examples.} than \eqref{eq:T_solution}, finding the left inverse $\mc{T}^*$ is very difficult both analytically and numerically; \cite{andrieu2021}.
\end{itemize}

To circumvent these challenges, it is reasonable to approximate these maps using neural networks. 

\begin{figure}[!ht]
    \begin{center}
        \begin{tikzpicture}
            \node (x) at (0,0) {$x^i(t_k)$};
            \node[fill=green!30!white!80!blue!70, rectangle, draw, minimum width=1.5cm, minimum height=1cm] (T1) at (2.5,1) {\Large $\hat{\mc{T}}_\theta$};
            \node (z) at (5,0) {$z^i(t_k)$};
            \node[fill=green!10!white!80!red!70, rectangle, draw, minimum width=1.5cm, minimum height=1cm] (T2) at (2.5,-1) {\Large $\hat{\mc{T}}_\eta^*$};
            
            \path[-latex, thick]
	            (x) edge (T1)
	            (T1) edge (z)
	            (z) edge (T2)
	            (T2) edge (x)
            ;
            \draw [very thick, decorate, decoration = {brace}] (1.75,1.65) --  (3.25,1.65);
            \node at (2.5,2) {\small $\frac{\partial \hat{\mc{T}}_\theta}{\partial x}(x^i(t_k)) f(x^i(t_k))=A\hat{\mc{T}}_\theta(x^i(t_k))+Bh(x^i(t_k))$};
            
        \end{tikzpicture}
    \end{center}
    \caption{Architecture for learning the transformation $\mc{T}$ and its inverse $\mc{T}^*$ using neural networks with parameters $\theta$ and $\eta$.}
    \label{fig:pinn_scheme}
\end{figure}
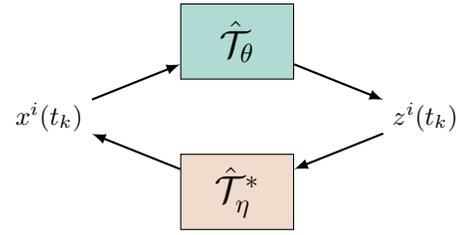

Let $\hat{\mc{T}}_{\theta}$ and $\hat{\mc{T}}^*_{\eta}$ be the parametrized neural networks that approximate $\mc{T}$ and $\mc{T}^*$, respectively. Here, $\theta, \eta$ are vectors containing all the weights and biases of each neural network, respectively, and can be considered as learning parameters for the nonlinear regression problem. In the following subsections, we describe our method, illustrated in Figure~\ref{fig:pinn_scheme}, for learning $\mc{T}$ and $\mc{T}^*$ through neural networks $\hat{\mc{T}}_{\theta}$ and $\hat{\mc{T}}^*_{\eta}$. 

\subsection{Generating Data for Training} \label{subsec_datagen}
Since the system trajectories for arbitrary initial conditions can be obtained numerically by solving the nonlinear system \eqref{eq:sys_x} for $x$ and the linear system \eqref{eq:sys_z} for $z$, one can pose the problem of learning $\theta$ and $\eta$ as a nonlinear regression over the simulated data trajectories on a finite time horizon $T>0$. The steps to generate these trajectories are described below:

\begin{enumerate}
    \item Define a set $\mc{X}^{\train} \times \mc{Z}^{\train}\subset \mc{X}\times \mc{Z}$ from which the initial conditions are chosen for training, where $\mc{Z}\subset\bb{R}^{n_z}$. For some $p\in\bb{N}$, choose a set of initial conditions 
    \[
        (x_0^1,z_0^1),\dots,(x_0^p,z_0^p) \in \mc{X}^{\train} \times \mc{Z}^{\train}.
    \]
    \item Simulate \eqref{eq:sys_x} and \eqref{eq:sys_z} with these initial conditions and generate sampled trajectories from $t_0=0$ to $t_{\tau-1}=T$
    \[
        x^i(t_k)\doteq x(t_k;x_0^i) ~\text{and}~ z^i(t_k)\doteq z(t_k;z_0^i)
    \]
    for $k=0,1,2,\dots,\tau-1$ and $i=1,\dots,p$.
        
    \item Partition the data samples into regression points $\ms{P}_r\subset\{0,\dots,\tau-1\}$ and physics points $\ms{P}_p\subset\{0,\dots,\tau-1\}$ such that $\ms{P}_r\cap \ms{P}_p=\emptyset$.
\end{enumerate} 

\begin{Remark} \label{remark:data_gen}
We provide the following guidelines for generating synthetic data trajectories:
\begin{enumerate}[(i)]
    \item The initial conditions $x_0^1,\dots,x_0^p$ can be chosen using the Latin hypercube sampling method; see \cite{ramos2020}.
    
    \item Choosing $z_0^1,\dots,z_0^p$ arbitrarily results in large regression errors for the initial time samples until the effect of the initial condition vanishes in $z(t;z_0^i)$ due to $A$ being Hurwitz. To avoid this, we follow a technique suggested by \cite{buisson2022}: 
    %
        (a)~Arbitrarily choose $p$ non-zero points $z_\tau^1,\dots,z_\tau^p$ in $\mc{Z}\subset\bb{R}^{n_z}$, where $\tau<0$ is such that 
        $
            \|\exp(A(t- \tau)) z_\tau^i\| \leq \epsilon
        $
        for some small $\epsilon>0$ and $t=0$. 
        Solving this inequality for $\tau$ gives
        \[ 
            \tau \leq \frac{1}{\lambda_{\min}(A)} \ln\left(\frac{\epsilon}{\cond(V) \bar{z}_\tau} \right)
        \]
        where $\bar{z}_\tau = \max_{i} \|z_\tau^i\|$ and $V$ is obtained from the eigendecomposition $A=V\Lambda V^{-1}$.
%
        (b)~Simulate \eqref{eq:sys_x} from $x_0^1,\dots,x_0^p$ in backward time and obtain output trajectories $h(x(t;x_0^i))$ for $t\in[\tau,0]$.
%
        (c)~Simulate \eqref{eq:sys_z} from $z_\tau^1,\dots,z_\tau^p$ in forward time and obtain $z(t;z_\tau^i)$ for $t\in[\tau,0]$. 
        %
        (d)~Choose $z_0^i = z(0;z_\tau^i)$, for $i=1,\dots,p$, which is approximately equal to $\mc{T}(x_0^i)$.
    
    \item A simple way to partition the data samples into regression points $\ms{P}_r$ and physics points $\ms{P}_p$ is to, for instance, choose even samples for $\ms{P}_r$ and odd samples for $\ms{P}_p$. \hfill $\diamond$
\end{enumerate}
\end{Remark}

\subsection{Defining the Empirical Loss Function} 
The regression problem minimizes a loss function that accounts for the deviation of the neural network's output with respect to the training data generated previously. To this end, we can exploit both $x^i(t_k)$ and $z^i(t_k)$ for learning $\mc{T}$ and $\mc{T}^*$ because both trajectories can be generated easily. The empirical loss function is defined as a \emph{mean squared error}
\begin{align}
\label{regress_cost}
    \mc{L}_{\theta,\eta}(X,Z) &\doteq \frac{1}{p} \sum_{i=1}^p \frac{1}{|\ms{P}_r|} \sum_{k\in \ms{P}_r} \big\|z^i(t_k)-\hat{\mc{T}}_\theta (x^i(t_k))\big\|^2 \notag \\
     & \qquad\qquad + \chi \big\|x^i(t_k)-\hat{\mc{T}}_\eta^* (\hat{\mc{T}}_\theta(x^i(t_k)))\big\|^2
\end{align}
where $\chi>0$ is a hyperparameter that not only weights the loss function properly but also discounts for different units of measurement of $x^i(t_k)$ and $z^i(t_k)$. Also, $X\in\bb{R}^{pn_x \times \tau}$ and $Z\in\bb{R}^{pn_z \times \tau}$ are defined as
\begin{align*}
    X &\doteq \left[\ba{cccc}
    x^1(t_0) & x^1(t_1) & \dots & x^1(t_{\tau-1}) \\
    x^2(t_0) & x^2(t_1) & \dots & x^2(t_{\tau-1}) \\
    \vdots & \vdots & \ddots & \vdots \\
    x^p(t_0) & x^p(t_1) & \dots & x^p(t_{\tau-1})
    \ea\right] \\
    Z &\doteq \left[\ba{cccc}
    z^1(t_0) & z^1(t_1) & \dots & z^1(t_{\tau-1}) \\
    z^2(t_0) & z^2(t_1) & \dots & z^2(t_{\tau-1}) \\
    \vdots & \vdots & \ddots & \vdots \\
    z^p(t_0) & z^p(t_1) & \dots & z^p(t_{\tau-1})
    \ea\right].
\end{align*}

\subsection{Enforcing the PDE Constraint}
An additional requirement of the learning problem is that $\hat{\mc{T}}_{\theta}$ must satisfy the PDE \eqref{eq:pde}
for every sample in $\mc{X}^\train$.
Evaluating \eqref{eq:pde} for all the physics points $\ms{P}_p$, we define the \emph{mean squared residual} of the PDE \eqref{eq:pde} over $\mc{X}^\train$ as
\begin{align}
\label{phys_res}
    \mc{N}_\theta(X) &\doteq \frac{1}{p} \sum_{i=1}^p \frac{1}{|\ms{P}_p|} \sum_{k\in \ms{P}_p} \big\|\frac{\partial \hat{\mc{T}}_\theta}{\partial x}(x^i(t_k)) f(x^i(t_k)) \notag\\ 
    & \qquad\qquad\qquad - A\hat{\mc{T}}_\theta(x^i(t_k)) - Bh(x^i(t_k))\big\|^2 
\end{align}
Enforcing the PDE constraint essentially avoids overfitting on the training samples and improves generalization by regularizing the neural network $\hat{\mc{T}}_\theta$.

\subsection{Supervised Physics-Informed Learning Problem}
By dedicating one part of the data for minimizing the mean squared error \eqref{regress_cost} and the other part for making the mean squared residual \eqref{phys_res} equal to zero, the \emph{supervised physics-informed learning problem} is formulated as:
\be \label{prob:PINN}
    \min_{\theta, \eta}  \mc{L}_{\theta,\eta}(X,Z) ~
    \text{subject to} ~ \mc{N}_\theta(X)=0.
\ee
Note that \eqref{prob:PINN} can be posed as
\be \label{prob:opt_final}
    \min_{\theta,\eta} \mc{L}_{\theta,\eta}(X,Z) + \lambda \mc{N}_\theta (X)
\ee
for a sufficiently large Lagrange multiplier $\lambda>0$ that discounts for the constraint $\mc{N}_\theta(X)=0$.

\subsection{Testing the Learned Model on a Different Dataset}

Once the neural networks $\hat{\mc{T}}_\theta$ and $\hat{\mc{T}}_\eta^*$ are trained, we evaluate the model's performance on the testing dataset $\mc{X}^{\test}\times \mc{Z}^{\test}\subset\mc{X}\times\mc{Z}$. It must be that the testing dataset is distinct from the training dataset for a fair evaluation of the performance. 
Moreover, we select multiple instances of testing dataset to tune the hyperparameters $\chi$ and $\lambda$ of the trained neural networks.
Among the two, $\lambda$ is a critical hyperparameter in \eqref{prob:opt_final} that largely impacts the satisfaction of the PDE constraint and, hence, the quality of the training.

\section{Evaluating the Performance of the Learned KKL Observer}
\label{sec_eval}

The neural networks $\hat{\mc{T}}_\theta$ and $\hat{\mc{T}}_\eta^*$ are mere approximations of $\mc{T}$ and $\mc{T}^*$, respectively. Thus, the performance of the observer will be influenced by the approximation error. Moreover, the model \eqref{eq:sys_x} of the state dynamics and sensors is never perfect in real-world applications, and there are several underlying uncertainties that could influence the state estimation. In this section, we provide robustness guarantees for the estimation error under both the approximation error and the system uncertainties. We also provide a metric to assess the generalization capability of the observer beyond the training data and discuss the specific features of the proposed learning method that avoid overfitting and enable better generalization as compared to other techniques.

\subsection{Robustness to the Approximation Error}

Given that the activation functions of the neural network are Lipschitz continuous, it can be shown that $\hat{\mc{T}}_\eta^*$ is also Lipschitz, i.e., there exists $\ell^*$ such that, for every $\hat{z},z\in\bb{R}^{n_z}$,
\be \label{eq:inverse_T_lipschitz}
    \|\hat{\mc{T}}_\eta^*(\hat{z}(t))-\hat{\mc{T}}_\eta^*(z(t))\| \leq \ell^* \|\hat{z}(t)-z(t)\|.
\ee
Specifically, we remark that ReLU networks are Lipschitz continuous, which is particularly important because we consider such a network in Section~\ref{sec_exp}. It is important to further remark that theoretical computation of the Lipschitz constant turns out to be quite conservative in practice. Although an NP-hard problem, empirically estimating a minimal Lipschitz constant of neural networks has been investigated extensively in the machine learning community \cite{szegedy2013, virmaux2018, fazlyab2019, jordan2020}.

For any $z\in\bb{R}^{n_z}$, $\mc{T}^*(z)$ can be written as
\be \label{eq:inverse_T}
\mc{T}^*(z) = \hat{\mc{T}}_\eta^*(z) + \mc{E}^*(z)
\ee 
where $\mc{E}^*(z)$ is the approximation error of $\hat{\mc{T}}^*$ at $z$.
Because the state space $\mc{X}\subset\bb{R}^{n_x}$ is bounded, $h(\cdot)$ is a smooth map, and $A$ is Hurwitz, there exists a compact set $\mc{Z}\subset\bb{R}^{n_z}$ containing the trajectory $z(t;\mc{T}(x_0))$ of \eqref{eq:sys_z} for every $t\geq 0$ and every $x_0\in\mc{X}$. Thus, as a consequence of \eqref{eq:uniformly_injective} and \eqref{eq:inverse_T_lipschitz}, there exists a finite approximation bound $\epsilon^*>0$ satisfying 
\be \label{eq:finite_approx_bound}
\epsilon^* = \sup_{z\in\mc{Z}} \|\mc{E}^*(z)\|.
\ee
There have been several attempts \cite{sokolic2017, kawaguchi2017, jakubovitz2019, cao2019} to estimate $\epsilon^*$ and to show that it can be reduced by improving the design and learning technique of the neural network, and also by increasing the size of the dataset $(X,Z)$ (see \cite{mohri2018}).

Using \eqref{eq:inverse_T}, we can write the KKL observer \eqref{eq:kkl_observer} as
\be \label{eq:kkl_observer2}
    \colsep=2pt\ba{ccl}
        \dot{\hat{z}} &=& A\hat{z} + By; \quad \hat{z}(0) = \hat{z}_0 \\
        \hat{x} &=& \hat{\mc{T}}_\eta^*(\hat{z}) + \mc{E}^*(\hat{z})
    \ea
\ee
where the approximation error $\mc{E}^*(\hat{z})$ is an unknown signal.

\begin{Proposition}
Subject to Assumptions~\ref{assumption_1} and \ref{assumption_2}, there exist positive constants $b,c>0$ such that the estimation error $\tilde{x}(t)=\hat{x}(t)-x(t)$ of \eqref{eq:kkl_observer2} satisfies
\be \label{eq:error_bound}
\|\tilde{x}(t)\| \leq b e^{-ct} + \epsilon^*, \quad\forall t\in\bb{R}_{\geq 0}
\ee
where
$
\epsilon^*
$
is given in \eqref{eq:finite_approx_bound}.
\end{Proposition}
\begin{proof}
We have
\begin{align} \label{eq:proof_ineq1}
    \|\tilde{x}(t)\| &= \|\hat{\mc{T}}_\eta^*(\hat{z}(t))-\mc{T}^*(z(t))\| \nonumber \\ 
    &= \|\hat{\mc{T}}_\eta^*(\hat{z}(t))-\hat{\mc{T}}_\eta^*(z(t))-\mc{E}^*(z(t))\| \nonumber \\ 
    &\leq  \|\hat{\mc{T}}_\eta^*(\hat{z}(t))-\hat{\mc{T}}_\eta^*(z(t))\| + \|\mc{E}^*(z(t))\| \nonumber \\ 
    &\leq  \ell^* \|\hat{z}(t) - z(t)\| + \epsilon^*
\end{align}
where the first step is due to \eqref{eq:inverse_T}, the second step is due to the triangle inequality, and the last step is due to \eqref{eq:inverse_T_lipschitz} and \eqref{eq:finite_approx_bound}. Since $A$ is Hurwitz, there exist $a,c>0$ such that
$
\|\hat{z}(t) - z(t) \| \leq ae^{-ct},
$
which completes the proof.
\end{proof}

\subsection{Robustness to Model Uncertainties and Sensor Noise}

Consider a nonlinear system
\be \label{eq:sys_x2}
\dot{x} = f(x) + w; \quad y=h(x) +v
\ee
where $w(t)\in\bb{R}^{n_x}$ and $v(t)\in\bb{R}^{n_y}$ are unknown but essentially bounded signals. In \eqref{eq:sys_x2}, the functions $f(\cdot)$ and $h(\cdot)$ represent the model of the system, and $w(t)$ represent model uncertainties and $v(t)$ the sensor noise.

We remark that the design method of KKL observers as presented in Sections~\ref{sec_KKLprelim} and \ref{sec_learning} remains the same for \eqref{eq:sys_x2}. However, to better attenuate the effects of uncertainties and noise, one can seek an $\mc{H}_\infty$-based design \cite{zemouche2016} of matrices $A$ and $B$ in the linear part of the KKL observer under the constraints that $A$ is Hurwitz and $(A,B)$ is controllable.

\begin{Proposition} \label{prop:robust_noise}
    Let Assumptions~\ref{assumption_1} and \ref{assumption_2} hold. Then, if $\|w\|_{L^\infty} \leq \bar{w}$ and $\|v\|_{L^\infty} \leq \bar{v}$, for every $t\in\bb{R}_{\geq 0}$, there exist positive constants $b,c,\alpha_1,\alpha_2>0$ such that the estimation error $\tilde{x}(t)=\hat{x}(t)-x(t)$ of \eqref{eq:kkl_observer2} satisfies
    \be \label{eq:error_bound_2}
        \|\tilde{x}(t)\| \leq b e^{-ct} + \alpha_1 \bar{w} + \alpha_2 \bar{v} + \epsilon^*, \quad\forall t\in\bb{R}_{\geq 0}
    \ee
    where
    $
    \epsilon^*
    $ is given in \eqref{eq:finite_approx_bound}.
\end{Proposition}
\begin{proof}[Proof idea]
    The proof follows from \eqref{eq:proof_ineq1} and the linear analysis of the error $\hat{z}(t)-z(t)$.
\end{proof}

Given that the model uncertainties and sensor noise are bounded, the above result shows that the KKL observer is robust in terms of input-to-state stability of the estimation error; see \cite{sontag1995}. Moreover, it can as well be shown that the constants $b,c,\alpha_1,\alpha_2$ in \eqref{eq:error_bound_2} are computable because of the linear dynamics of the KKL observer.

\subsection{Assessing the Observer's Generalization Capability} \label{subsec_gen_metric}

Another key contribution in this paper is to evaluate the performance of the learned KKL observer even when the true initial condition of the system in real-time is far from the training region $\mc{X}^{\train}$.
To this end, we define a metric quantifying the generalization capability of the trained model in Figure~\ref{fig:pinn_scheme} for the KKL observer. This metric compares the estimation errors resulting from the training and the testing phases, and describes how the error varies as a function of the distance between the two sets $\mc{X}^{\train}$ and $\mc{X}^{\test}$. 

Let the testing region $\mc{X}^{\test} \subset \mc{X}\setminus \mc{X}^{\train}$,
and consider a set of points
$
\{\xi_0^j : j = 1, \dots, q\} \in \mc{X}^{\test}
$
that, for every $j\in\{1,\dots,q\}$, satisfy
$
d(\xi_0^j,\mc{X}^\train) = \delta
$,
for some $\delta>0$, where 
\[
d(\xi_0^j,\mc{X}^\train)\doteq \inf\limits_{x_0\in\mc{X}^\train} \| x_0 - \xi_0^j \|.
\]
The \textit{empirical generalization error} $\mathrm{G}_{\text{emp}} (\delta)$ is defined as
\begin{align} \label{eq:gen_metric}
    \mathrm{G}_{\text{emp}}(\delta)\doteq |\mathrm{E}_{\test}(\delta) - \mathrm{E}_{\train}|
\end{align}
where
\begin{align*}
    \mathrm{E}_{\test}(\delta) &\doteq  \frac{1}{q} \sum_{j=1}^q \frac{1}{\tau}\sum\limits_{k=0}^{\tau-1} \frac{\| \hat{x}(t_k; \hat{\xi}_0^j) - x(t_k; \xi_0^j(\delta)) \|^2}{\|x(t_k; \xi_0^j(\delta))\|^2} \\
    \mathrm{E}_{\train} &\doteq \frac{1}{p}\sum\limits_{i=1}^{p}\frac{1}{\tau}\sum\limits_{k=0}^{\tau-1} \frac{\| \hat{x}(t_k; \hat{x}_0^i) - x(t_k; x_0^i) \|^2}{\|x(t_k;x_0^i)\|^2}
\end{align*}
with $\hat{\xi}_0^j$ and $\hat{x}_0^i$ chosen sufficiently close to $\xi_0^j(\delta)$ and $x_0^i$, respectively, to avoid the errors accumulated in the observer's transient.
Notice that $\mathrm{E}_\test$ denotes the normalized mean estimation error variance of multiple test trajectories initialized at $\delta$-distance from $\mc{X}^\train$, whereas $\mathrm{E}_\train$ denotes the normalized mean estimation error variance of all the training trajectories.

In short, during the testing phase, we select $q$ initial points $\xi_0^j$ that are $\delta$-distant from the training region, where $\delta\in\{\delta_1,\dots,\delta_m\}$ with $0<\delta_1<\dots<\delta_m$. Then, for each $\delta_i$, the change in the normalized testing error variance $\mathrm{E}_{\test}(\delta_i)$ provides an empirical quality measure \eqref{eq:gen_metric} on the generalization capability of the learned KKL observer. 

\subsection{Discussion on the Observer's Generalization Capability}

Since Assumptions~\ref{assumption_1} and \ref{assumption_2} ensure uniform injectivity of $\mc{T}$, and $\mc{T}$ satisfies the PDE \eqref{eq:pde}, the inverse $\mc{T}^*$ exists and is unique. Thus, the data samples used in the training are of the form $(x,\mc{T}(x))$ and $(z,\mc{T}^*(z))$, which entails that the problem \eqref{prob:opt_final} is a \textit{realizable} learning task that is \textit{probably approximately correct} (PAC) learnable \cite{shalev2014}. Then, one of the sources for non-zero generalization error is the fact that the training data $(X,Z)$ induced loss $\mc{L}_{\theta,\eta}(X,Z)$ in \eqref{regress_cost} is an approximation of the actual loss 
\[\ba{r}
\displaystyle \bar{\mc{L}}_{\theta,\eta}(x,z) \doteq \int_{\mc{X}} \int_0^T \|z(t;\mc{T}(\xi))-\mc{T}(x(t;\xi))\| \qquad \\ [1em]
+ \chi \|x(t;\xi) - \mc{T}^*(\mc{T}(x(t;\xi)))\| \text{d}t \, \text{d}\mu(\xi)
\ea\]
where $\mu$ is a measure on $\mc{X}$.

In our formulation, an unlimited amount of synthetic data can be generated using the method described in Section~\ref{subsec_datagen}, which enables one to enhance the generalization capability of the learned KKL observer and improve its performance.
However, it is not practical to utilize arbitrarily large amount of data for training. Therefore, under the same training data size, a key feature that makes the supervised PINN to have better generalization capability than the neural network architectures of \cite{ramos2020,peralez2021,buisson2022} is the regularization with the PDE \eqref{eq:pde}, which reduces the search space of the hypothesis and avoids overfitting on the training data.

In the unsupervised AE architecture of \cite{buisson2022}, the neural network is also regularized by the PDE \eqref{eq:pde}. However, unlike \eqref{regress_cost}, the loss function of \cite{buisson2022} doesn't include additional regression term that accounts for the deviation between $z^{i}(t_k)$ and $\hat{\mc{T}}_\theta (x^i(t_k))$. This is very important because without the explicit supervision to connect the system's state space $\mc{X}$ to the observer's state space $\mc{Z}$, the AE will minimize the reconstruction loss on a limited number of training samples $x^i(t_k)$, which may belong to a larger hypothesis space. Thus, the unsupervised AE of \cite{buisson2022} makes the neural network overfit upon the partial training data, i.e., only in the $x$-domain, and hinders the generalization on the unseen data. 
In the extreme case, without the PDE regularization, if the decoder is complex enough, one could essentially recover the $x$ sample even from noise, and the learned left inverse $\hat{\mc{T}}_\eta^*$ can as well be arbitrary \cite{vandenoord2016}.

\begin{figure*}[h!t]
        \centering
        \begin{subfigure}[t]{0.32\textwidth}
        \includegraphics[width=\textwidth]{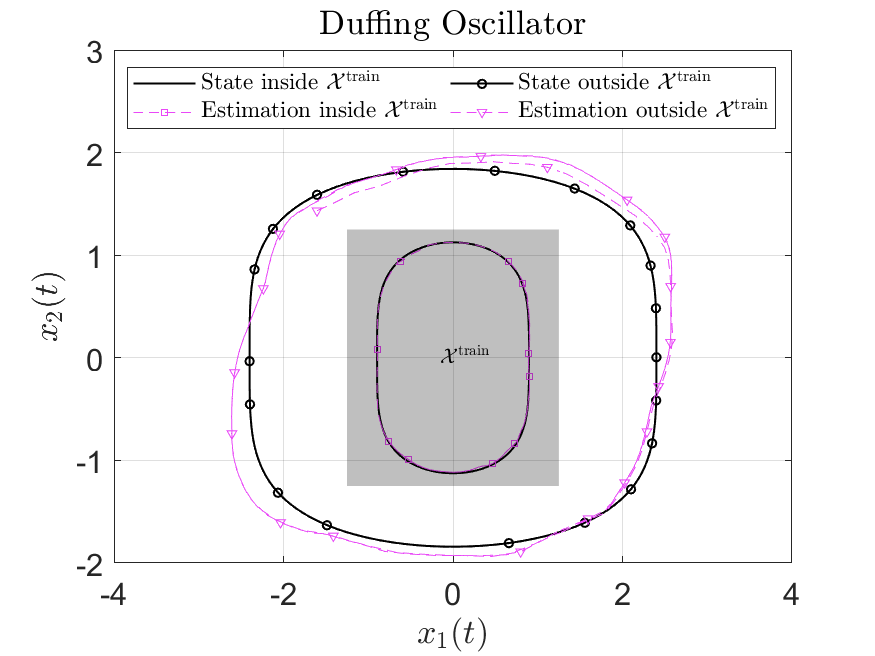}
        
        \includegraphics[width=\textwidth]{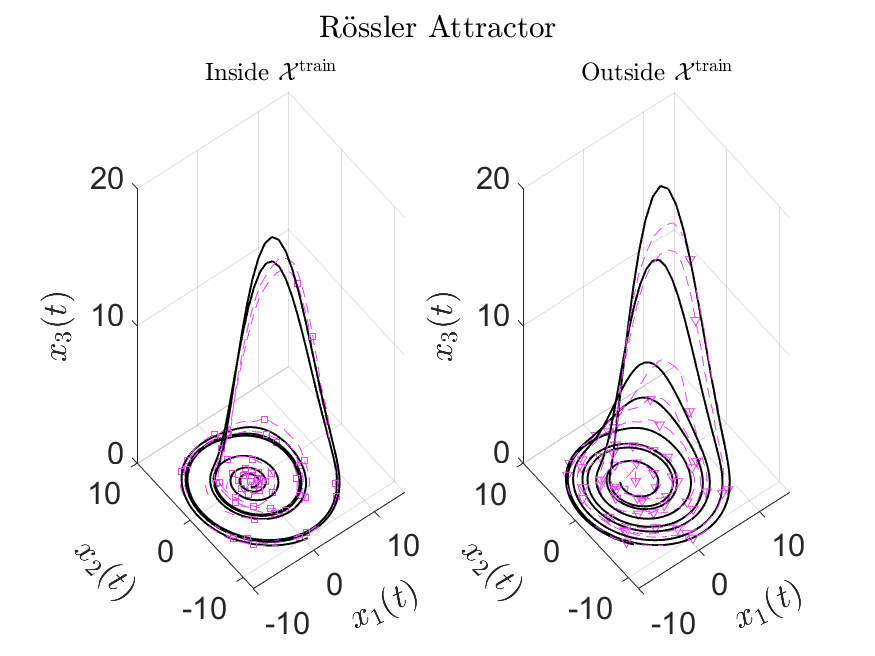}
        \caption{}
        \label{fig:in_and_out_train}
        \end{subfigure}
        \begin{subfigure}[t]{0.32\textwidth}
        \includegraphics[width=\textwidth]{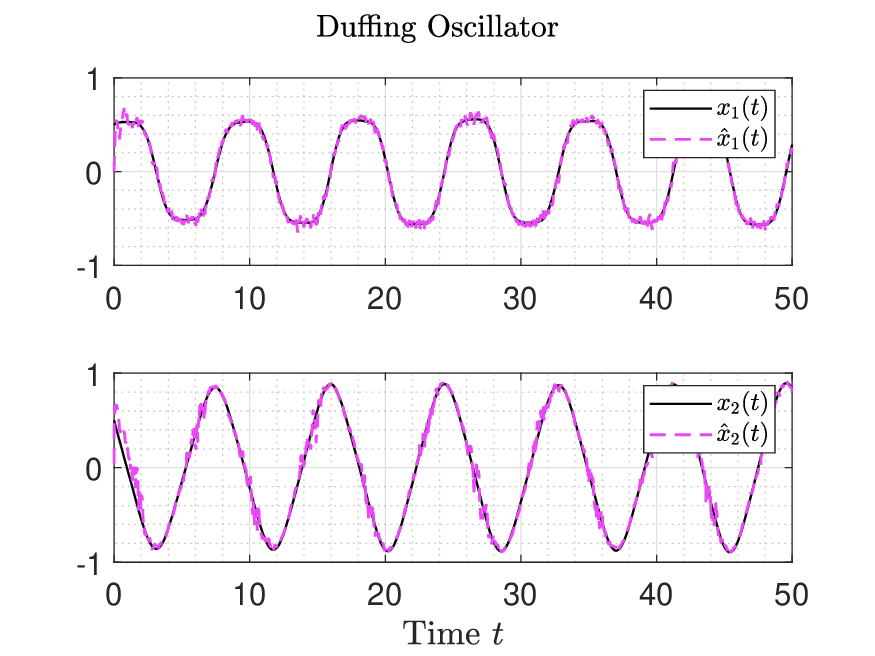}
        
        \includegraphics[width=\textwidth]{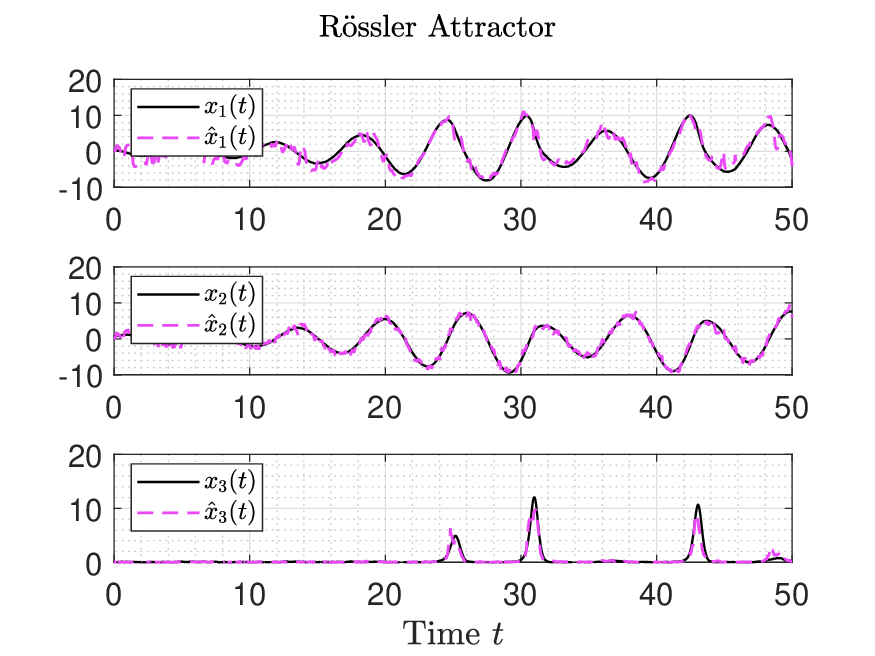}
        \caption{}
        \label{fig:estimation_noise}
        \end{subfigure}
        \begin{subfigure}[t]{0.32\textwidth}
        \includegraphics[width=\textwidth]{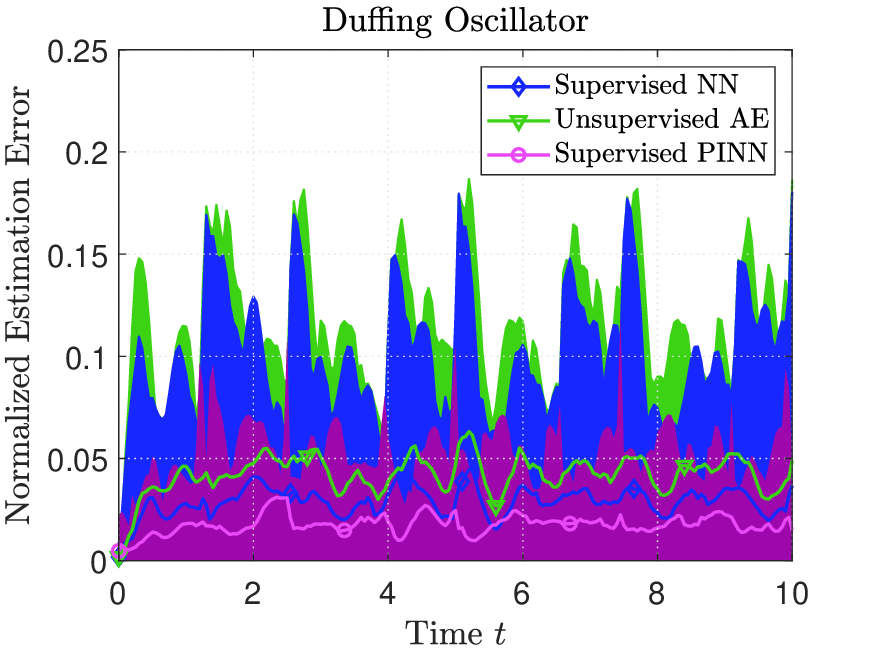}  
        
        \includegraphics[width=\textwidth]{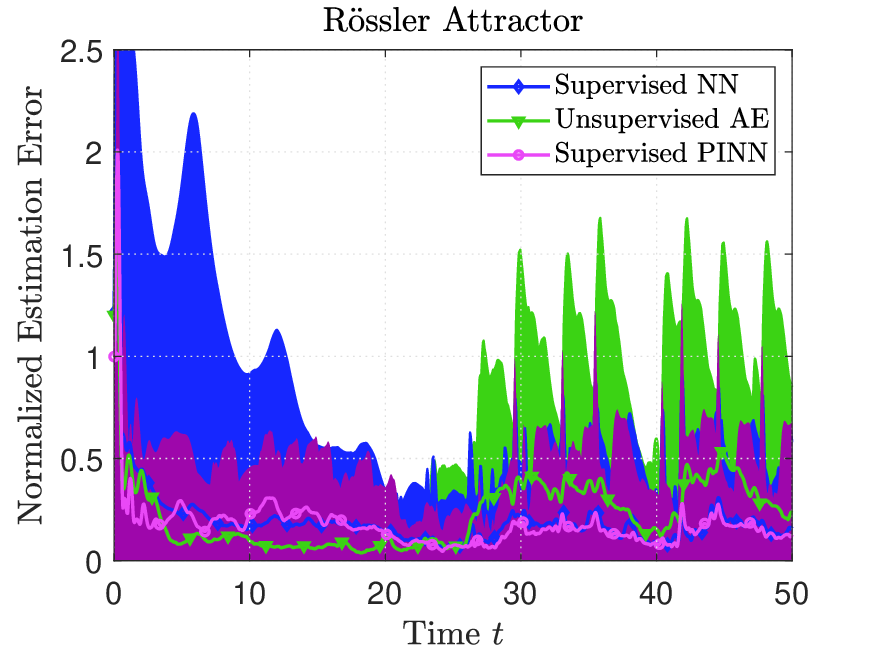}  
        \caption{}
        \label{fig:error_range_test}
        \end{subfigure}
        \caption{(a)~Phase portrait of the estimation performance when the system is initialized inside and outside the training region. (b)~Estimation performance in the presence of model uncertainties and sensor noise. (c)~Comparison of our method with others in terms of the range of normalized estimation errors and their averages for 50 state trajectories initialized outside the training region.}
\end{figure*}

\section{Experimentation and Testing} \label{sec_exp}

Performance of the proposed supervised PINN-based KKL observer is numerically tested under different scenarios. First, we test its performance under approximation errors when the state trajectory is initialized outside the training region $\mc{X}^\train$. Second, we test its performance under model uncertainties and sensor noise and demonstrate the robustness of the proposed observer. Third, we examine the estimation error trajectories for multiple experiments where the system's state is always initialized randomly outside $\mc{X}^\train$. We show that the proposed supervised PINN-based KKL observer demonstrates better performance than 1)~supervised NN \cite{ramos2020, peralez2021, buisson2022} and 2)~unsupervised AE \cite{buisson2022}. Finally, we compare the empirical generalization error resulting from all these techniques and demonstrate that our method exhibits better generalization capabilities. 

For the experimentation and testing, we consider the following nonlinear oscillators:
\begin{itemize}
    \item {Reverse Duffing oscillator}
    \be \label{eq:revduffing}
        \dot{x}_1 = x_2^3, \quad \dot{x}_2 = -x_1; \quad
        y = x_1.
    \ee
    \item {R\"{o}ssler attractor}
    \be \label{eq:rossler}
    \colsep=2pt\ba{rclrcl}
        \dot{x}_1 &=& -x_2 -x_3, \quad
        &\dot{x}_2 &=& x_1 + a x_2 \\
        \dot{x}_3 &=& b + x_3(x_1 - c); \quad 
        &y &=& x_2
    \ea
    \ee
    where the parameters $a=0.2$, $b=0.2$, and $c=5.7$. 
\end{itemize}

\subsection{Experimental Setup for Training and Testing}
For both \eqref{eq:revduffing} and \eqref{eq:rossler}, we follow the data generation and sampling procedure described in Section~\ref{subsec_datagen}. For reverse Duffing oscillator, $\mathcal{X}^{\train}=[-1,1]^2$. For R\"{o}ssler attractor, $\mathcal{X}^{\train}=[-1,1]^3$. We generate $\{x_0^1,...,x_0^p\}$ using Latin hypercube sampling method. The initial conditions $\{z_0^1,\dots,z_0^p\}$ are generated using Remark~\ref{remark:data_gen}(ii). 
Runge-Kutta-4 is used as the numerical ODE-solver for \eqref{eq:revduffing}-\eqref{eq:rossler} over a time horizon $[0,50]$. 

The architecture of both neural networks $\hat{\mc{T}_\theta}$ and $\hat{\mc{T}^*_\eta}$ in Figure \ref{fig:pinn_scheme} is chosen to be a multi-layer perceptron with five hidden layers, where each layer has 50 neurons with ReLU activation function. We use normalization and denormalization layer for data standardization in order to facilitate the training.
Training is further facilitated by a learning rate scheduler. All models in this section are trained using the \texttt{Adam} optimization algorithm with a batch size of 32.
In the testing stage, initial conditions are generated outside the training domain, from which \eqref{eq:revduffing} and \eqref{eq:rossler} are then simulated.  For the code and other details, please refer to our repository\footnote{\url{https://github.com/Mudhdhoo/ACC_KKL_Observer}}.

The matrices of the KKL observer are chosen as follows:
\[
A = -\diag(1,2,\dots,n_z), \quad B=1_{n_z}
\]
where $n_z=n_y(2n_x+1)$, $\diag()$ denotes a diagonal matrix, and $1_{n_z}$ is a vector of ones with dimensions $n_z\times 1$. Notice that $n_y=1$ for both \eqref{eq:revduffing} and \eqref{eq:rossler}.

\subsection{Experimental Results}

In the following, we present several experimental results and compare our method \textit{supervised PINN} with \textit{supervised NN} \cite{ramos2020,peralez2021,buisson2022} and \textit{unsupervised AE} \cite{buisson2022}.

\subsubsection{Testing the supervised PINN-based KKL observer outside the training region}
We train the supervised PINN inside the training regions for both \eqref{eq:revduffing} and \eqref{eq:rossler}. We test it outside the training region. Figure~\ref{fig:in_and_out_train} demonstrates the estimation performance of the learned KKL observer when the true system is initialized inside the training region $\mc{X}^\train$ and outside the training region. Despite an expected deterioration of the state estimation outside the training region, the observer's performance is satisfactory as it is able to follow the true state with a small error. 

\subsubsection{Testing the supervised PINN-based KKL observer under model uncertainties and sensor noise}

We randomly initialize the state trajectories inside the training region $\mc{X}^\train$, where the initial point is different from the initial points in the training dataset $X$. We consider $w(t)\sim \mc{N}(0, 0.1)$ and $v(t)\sim \mc{N}(0,0.1)$ for \eqref{eq:revduffing}, $w(t)\sim \mc{N}(0,1)$ and $v(t)\sim \mc{N}(0,1)$ for \eqref{eq:rossler}. Figure~\ref{fig:estimation_noise} shows the true and estimated state trajectories, and demonstrates that the learned KKL observer is stable under uncertainties and noise as stated in Proposition~\ref{prop:robust_noise}.

\subsubsection{Estimation errors for multiple state trajectories initialized outside the training region}
We initialize the systems \eqref{eq:revduffing} and \eqref{eq:rossler} at 50 points that are randomly generated outside the training region. We run the KKL observers that are learned according to supervised NN, unsupervised AE, and our method supervised PINN. To show the merits of each learning scheme, we compare \textit{normalized estimation error} trajectories
\[
e_i(t) = \frac{\|\hat{x}^i(t) - x^i(t)\|}{\|x^i(t)\|}; \quad i = 1, \ldots, 50.
\]
Figure~\ref{fig:error_range_test} demonstrates the error ranges and the average ($\sum_{i=1}^{50} e_i(t)/50$) for each learning scheme. For the reverse Duffing oscillator, our method yields lowest maximum and average error for all times. For the R\"{o}ssler attractor, the overall performance of our method is better than both the supervised NN and unsupervised AE. The supervised NN performs worse in the beginning, which is before the bifurcation of the R\"{o}ssler attractor, because it fails to capture some trajectories that are initialized outside the training region. On the other hand, the unsupervised AE performs worse after the bifurcation because it is not very sensitive to changes in the $z$-domain that correspond to the bifurcation in the $x$-domain.

\subsubsection{Comparison of the empirical generalization error for multiple learning schemes}

\begin{figure}[!]
        \centering
        \includegraphics[width=0.45\textwidth]{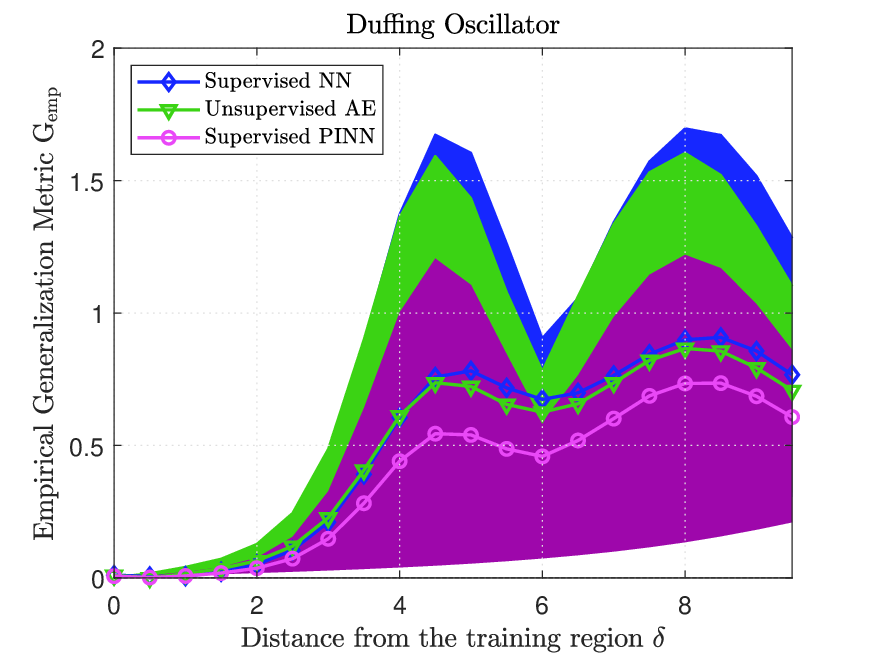}
        \caption{Comparison of the empirical generalization error as the initial state of reverse Duffing oscillator is at a distance $\delta$ from the training region $\mc{X}^\train$.}
        \label{fig:gen_metric}
\end{figure}

We choose multiple initial points outside the training region for each $\delta_i>0$ in the testing phase as described in Section~\ref{subsec_gen_metric}. We only consider reverse Duffing oscillator \eqref{eq:revduffing} for this experiment. We choose multiple $\delta_i\in\{0.5,1,1.5,\dots,10\}$, and, for each $\delta_i$, we choose 10 initial points in circular formation centered around $[-1,1]^2$ outside $\mc{X}^\train$. 
Figure~\ref{fig:gen_metric} illustrates the comparison of different learning schemes in terms of empirical generalization error. For all $\delta_i$, it can be seen that supervised PINN yields smaller generalization errors.

\section{Discussion and Future Outlook} \label{sec_conc}

We proposed a novel supervised physics-informed learning method to design Luenberger or KKL observers for autonomous nonlinear systems. The proposed method learns the nonlinear transformation map required to transform the system to the observer's coordinates and satisfies a certain PDE constraint. Additionally, the inverse of the transformation map is learned to obtain the state estimate in the original state space. To learn both the transformation map and its inverse, we trained a physics-informed neural network architecture on synthetic data generated by numerically solving both the system and the observer. The PDE constraint acts as a physical invariant that regularizes the neural network, reducing the hypothesis's search space. We demonstrated that the KKL observer designed with our method is robust to neural network's approximation error, model uncertainties, and sensor noise. The proposed method also exhibits better generalization properties than other methods due to the PDE regularization and the regression loss in the observer's coordinates. We validated our results on reverse Duffing oscillator and R\"{o}ssler attractor.

While we discussed the generalization capability of the proposed learning-based observer design method in detail, theoretical guarantees on its generalizability remain an open problem. Additionally, designing a KKL observer optimally to improve its robustness to model uncertainties and sensor noise is left for future work. We also recognize the potential of alternative methods such as operator learning \cite{kissas2022} to learn the non-linear transformation map, which is a prospect to be explored. Furthermore, the proposed method can be extended beyond KKL observers to obtain the triangular form of nonlinear systems required in designing high-gain and backstepping observers. 
In conclusion, the proposed learning-based observer design method can be a promising solution to address the challenging problem of designing observers for nonlinear systems.

\bibliographystyle{IEEEtran}
\bibliography{biblio}

\begin{thebibliography}{10}
\providecommand{\url}[1]{#1}
\csname url@samestyle\endcsname
\providecommand{\newblock}{\relax}
\providecommand{\bibinfo}[2]{#2}
\providecommand{\BIBentrySTDinterwordspacing}{\spaceskip=0pt\relax}
\providecommand{\BIBentryALTinterwordstretchfactor}{4}
\providecommand{\BIBentryALTinterwordspacing}{\spaceskip=\fontdimen2\font plus
\BIBentryALTinterwordstretchfactor\fontdimen3\font minus
  \fontdimen4\font\relax}
\providecommand{\BIBforeignlanguage}[2]{{%
\expandafter\ifx\csname l@#1\endcsname\relax
\typeout{** WARNING: IEEEtran.bst: No hyphenation pattern has been}%
\typeout{** loaded for the language `#1'. Using the pattern for}%
\typeout{** the default language instead.}%
\else
\language=\csname l@#1\endcsname
\fi
#2}}
\providecommand{\BIBdecl}{\relax}
\BIBdecl

\bibitem{luenberger1964}
D.~G. Luenberger, ``Observing the state of a linear system,'' \emph{IEEE
  Transactions on Military Electronics}, vol.~8, no.~2, pp. 74--80, 1964.

\bibitem{shoshitaishvili1990}
A.~Shoshitaishvili, ``Singularities for projections of integral manifolds with
  applications to control and observation problems,'' in \emph{Theory of
  Singularities and its Applications}.\hskip 1em plus 0.5em minus 0.4em\relax
  American Mathematical Society, 1990, pp. 295--333.

\bibitem{shoshitaishvili1992}
------, ``On control branching systems with degenerate linearization,'' in
  \emph{IFAC Symposium on Nonlinear Control Systems}, 1992, pp. 495--500.

\bibitem{kazantzis1998}
N.~Kazantzis and C.~Kravaris, ``Nonlinear observer design using lyapunov’s
  auxiliary theorem,'' \emph{Systems \& Control Letters}, vol.~34, no.~5, pp.
  241--247, 1998.

\bibitem{krener2002b}
A.~J. Krener and M.~Xiao, ``Nonlinear observer design in the siegel domain,''
  \emph{SIAM Journal on Control and Optimization}, vol.~41, no.~3, pp.
  932--953, 2002.

\bibitem{kreisselmeier2003}
G.~Kreisselmeier and R.~Engel, ``Nonlinear observers for autonomous {Lipschitz}
  continuous systems,'' \emph{IEEE Transactions on Automatic Control}, vol.~48,
  no.~3, pp. 451--464, 2003.

\bibitem{andrieu2006}
V.~Andrieu and L.~Praly, ``On the existence of a
  {Kazantzis--Kravaris/Luenberger} observer,'' \emph{SIAM Journal on Control
  and Optimization}, vol.~45, no.~2, pp. 432--456, 2006.

\bibitem{andrieu2014}
V.~Andrieu, ``Convergence speed of nonlinear {Luenberger} observers,''
  \emph{SIAM Journal on Control and Optimization}, vol.~52, no.~5, pp.
  2831--2856, 2014.

\bibitem{engel2007}
R.~Engel, ``Nonlinear observers for {Lipschitz} continuous systems with
  inputs,'' \emph{International Journal of Control}, vol.~80, no.~4, pp.
  495--508, 2007.

\bibitem{bernard2017}
P.~Bernard, ``Luenberger observers for nonlinear controlled systems,'' in
  \emph{2017 IEEE 56th Annual Conference on Decision and Control (CDC)}.\hskip
  1em plus 0.5em minus 0.4em\relax IEEE, 2017, pp. 3676--3681.

\bibitem{bernard2018}
P.~Bernard and V.~Andrieu, ``Luenberger observers for nonautonomous nonlinear
  systems,'' \emph{IEEE Transactions on Automatic Control}, vol.~64, no.~1, pp.
  270--281, 2018.

\bibitem{andrieu2021}
V.~Andrieu and P.~Bernard, ``Remarks about the numerical inversion of injective
  nonlinear maps,'' in \emph{2021 60th IEEE Conference on Decision and Control
  (CDC)}, 2021, pp. 5428--5434.

\bibitem{bernard2022}
P.~Bernard, V.~Andrieu, and D.~Astolfi, ``Observer design for continuous-time
  dynamical systems,'' \emph{Annual Reviews in Control}, 2022.

\bibitem{ramos2020}
L.~d.~C. Ramos, F.~Di~Meglio, V.~Morgenthaler, L.~F.~F. da~Silva, and
  P.~Bernard, ``Numerical design of {Luenberger} observers for nonlinear
  systems,'' in \emph{2020 59th IEEE Conference on Decision and Control (CDC)},
  2020, pp. 5435--5442.

\bibitem{peralez2021}
J.~Peralez and M.~Nadri, ``Deep learning-based {Luenberger} observer design for
  discrete-time nonlinear systems,'' in \emph{2021 60th IEEE Conference on
  Decision and Control (CDC)}, 2021, pp. 4370--4375.

\bibitem{buisson2022}
M.~Buisson-Fenet, L.~Bahr, and F.~Di~Meglio, ``Towards gain tuning for
  numerical {KKL} observers,'' \emph{arXiv preprint arXiv:2204.00318}, 2022.

\bibitem{raissi2019}
M.~Raissi, P.~Perdikaris, and G.~E. Karniadakis, ``Physics-informed neural
  networks: A deep learning framework for solving forward and inverse problems
  involving nonlinear partial differential equations,'' \emph{Journal of
  Computational Physics}, vol. 378, pp. 686--707, 2019.

\bibitem{karniadakis2021}
G.~E. Karniadakis, I.~G. Kevrekidis, L.~Lu, P.~Perdikaris, S.~Wang, and
  L.~Yang, ``Physics-informed machine learning,'' \emph{Nature Reviews
  Physics}, vol.~3, no.~6, pp. 422--440, 2021.

\bibitem{sontag1995}
E.~D. Sontag and Y.~Wang, ``On characterizations of the input-to-state
  stability property,'' \emph{Systems \& Control Letters}, vol.~24, no.~5, pp.
  351--359, 1995.

\bibitem{brivadis2022}
L.~Brivadis, V.~Andrieu, P.~Bernard, and U.~Serres, ``Further remarks on {KKL}
  observers,'' \emph{HAL preprint HAL-03695863}, 2022.

\bibitem{szegedy2013}
C.~Szegedy, W.~Zaremba, I.~Sutskever, J.~Bruna, D.~Erhan, I.~Goodfellow, and
  R.~Fergus, ``Intriguing properties of neural networks,'' \emph{arXiv preprint
  arXiv:1312.6199}, 2013.

\bibitem{virmaux2018}
A.~Virmaux and K.~Scaman, ``Lipschitz regularity of deep neural networks:
  Analysis and efficient estimation,'' \emph{Advances in Neural Information
  Processing Systems}, vol.~31, 2018.

\bibitem{fazlyab2019}
M.~Fazlyab, A.~Robey, H.~Hassani, M.~Morari, and G.~Pappas, ``Efficient and
  accurate estimation of {Lipschitz} constants for deep neural networks,''
  \emph{Advances in Neural Information Processing Systems}, vol.~32, 2019.

\bibitem{jordan2020}
M.~Jordan and A.~G. Dimakis, ``Exactly computing the local {Lipschitz} constant
  of {ReLU} networks,'' \emph{Advances in Neural Information Processing
  Systems}, vol.~33, pp. 7344--7353, 2020.

\bibitem{sokolic2017}
J.~Sokoli{\'c}, R.~Giryes, G.~Sapiro, and M.~R. Rodrigues, ``Robust large
  margin deep neural networks,'' \emph{IEEE Transactions on Signal Processing},
  vol.~65, no.~16, pp. 4265--4280, 2017.

\bibitem{kawaguchi2017}
K.~Kawaguchi, L.~P. Kaelbling, and Y.~Bengio, ``Generalization in deep
  learning,'' \emph{arXiv preprint arXiv:1710.05468}, 2017.

\bibitem{jakubovitz2019}
D.~Jakubovitz, R.~Giryes, and M.~R. Rodrigues, ``Generalization error in deep
  learning,'' in \emph{Compressed sensing and its applications}.\hskip 1em plus
  0.5em minus 0.4em\relax Springer, 2019, pp. 153--193.

\bibitem{cao2019}
Y.~Cao and Q.~Gu, ``Generalization bounds of stochastic gradient descent for
  wide and deep neural networks,'' \emph{Advances in Neural Information
  Processing Systems}, vol.~32, 2019.

\bibitem{mohri2018}
M.~Mohri, A.~Rostamizadeh, and A.~Talwalkar, \emph{Foundations of Machine
  Learning}.\hskip 1em plus 0.5em minus 0.4em\relax MIT press, 2018.

\bibitem{zemouche2016}
A.~Zemouche, R.~Rajamani, B.~Boulkroune, H.~Rafaralahy, and M.~Zasadzinski,
  ``$\mathcal{H}_\infty$ circle criterion observer design for {Lipschitz}
  nonlinear systems with enhanced {LMI} conditions,'' in \emph{2016 American
  Control Conference (ACC)}, 2016, pp. 131--136.

\bibitem{shalev2014}
S.~Shalev-Shwartz and S.~Ben-David, \emph{Understanding Machine Learning: From
  Theory to Algorithms}.\hskip 1em plus 0.5em minus 0.4em\relax Cambridge
  University Press, 2014.

\bibitem{vandenoord2016}
A.~Van~den Oord, N.~Kalchbrenner, L.~Espeholt, O.~Vinyals, A.~Graves
  \emph{et~al.}, ``Conditional image generation with {PixelCNN} decoders,''
  \emph{Advances in Neural Information Processing Systems}, vol.~29, 2016.

\bibitem{kissas2022}
G.~Kissas, J.~H. Seidman, L.~F. Guilhoto, V.~M. Preciado, G.~J. Pappas, and
  P.~Perdikaris, ``Learning operators with coupled attention,'' \emph{Journal
  of Machine Learning Research}, vol.~23, no. 215, pp. 1--63, 2022.

\end{thebibliography}

\end{document}